\documentclass[12pt,reqno]{amsart}

\usepackage[foot]{amsaddr}

\makeatletter
\renewcommand{\email}[2][]{%
  \ifx\emails\@empty\relax\else{\g@addto@macro\emails{,\space}}\fi%
  \@ifnotempty{#1}{\g@addto@macro\emails{\textrm{(#1)}\space}}%
  \g@addto@macro\emails{#2}%
}
\makeatother
\title{The title}
\usepackage{graphicx}
\usepackage{xfrac}
\usepackage{amsfonts,amsmath,amssymb,epsfig,mathrsfs,bm}
\usepackage{ifthen,color,mathtools}
\usepackage{epic,eepic}
\usepackage{esint,stmaryrd}
\usepackage{a4wide,accents}
\usepackage{mathscinet}
\usepackage{cite}
\usepackage{comment}

\usepackage[top=3cm,bottom=2.5cm,left=2cm,right=2cm]{geometry}

\ProvideTextCommandDefault{\cprime}{\tprime}

\newtheorem{theorem}{\bf Theorem}[section]
\newtheorem{proposition}{\bf Proposition}[section]

\newtheorem{lemma}{\bf Lemma}[section]

\newtheorem{remark}{\bf Remark}[section]

\numberwithin{equation}{section}
\newcommand{\diver}{\operatorname{div}}
\newcommand{\dx}{\,{\rm d} {x}}
\newcommand{\vr}{\varrho}
\newcommand{\bProof}{{\bf Proof. }}
\newcommand{\intS}[1]{\int_{B_\vr^{+}} #1 \ \dx}
\newcommand{\intSS}[1]{\int_{B_{2\vr}^+} #1 \ \dx}

\newcommand{\intB}[1]{\int_{(\partial B_\vr)^{+}} #1 \ \d\sigma_x}
\newcommand{\intBp}[1]{\int_{\partial B^{+}_\vr} #1 \ \d\sigma_x}
\newcommand{\R}{{\mathbb R}}
\newcommand{\Id}{\mathrm{Id}}

\renewcommand{\d}{\,{\rm d}}
\newcommand{\pB}[1]{\partial B_{#1}}
\newcommand{\Gradp}{\nabla_p}

\begin{document}
\title{On the nonlinear thin obstacle problem}
\author[A.~Abbatiello]{Anna Abbatiello}
\address{Universit\`{a} degli Studi della Campania ``L. Vanvitelli'', Dipartimento di Matematica e Fisica, Viale A. Lincoln 5, 81100 Caserta, Italy}
\email{anna.abbatiello@unicampania.it}

\author[G.~Andreucci]{Giovanna Andreucci}
\author[E.~Spadaro]{Emanuele Spadaro}
\address{Sapienza University of Rome, Dipartimento di Matematica ``G.Castelnuovo", Piazzale Aldo Moro 5, 00185 Rome, Italy.}
\email{giovanna.andreucci@uniroma1.it, spadaro@mat.uniroma1.it}

\begin{abstract}
The thin obstacle problem or $n$-dimensional Signorini problem is a classical variational problem arising in several applications, starting with its first introduction in elasticity theory. The vast literature concerns mostly quadratic energies, whereas only partial results have been proved in the nonlinear case. In this paper we consider the thin boundary obstacle problem for a general class of nonlineraities and we prove the optimal $C^{1, \frac{1}{2}}$-regularity of the solutions in any space dimension. 
\end{abstract}

\maketitle
\section{Introduction}
We are interested in the boundary (thin) obstacle problems for a class of nonlinear functionals of the type:
\begin{equation}\label{min}
 \min_{u\in \mathcal{A}_g} \; \int_{B_{1}^+} f\left(\nabla u\right) \d x ,
\end{equation}
where, for any subset $E\subset \R^n$, we use the notation
\[
E^+= E\cap \{x_n>0\} \quad \textup{and}\quad
E':= E\cap \{x_n=0\},
\]
$B_1$ is the unit ball in $\R^n$ with $n\geq 2$ and
\begin{itemize}
\item the class of competitor functions is 
$$\mathcal{A}_g:=\Big\{u\in W_0^{1,\infty}{(B_1^+)}: u_{|B'_1}\geq 0, \; u_{|(\partial B_1)^+} = g_{|(\partial B_1)^+} 
\Big\},$$
with $g\in C^2(\R^n)$ prescribes the boundary values and satisfies $g_{|\R^{n-1}\times\{0\}}\geq 0$;

\item the nonlinear (non-quadratic) energy density $f:\R^n\to \R$ is of the form
\begin{equation}\label{formula-f}
f(p)=h(|p|)\quad 
\forall\; p \in \mathbb{R}^{n} ,
\end{equation}
with $h\in C^2(\R)$ satisfying
\begin{equation}\label{e.h''}
h(0)=h'(0)=0, \quad  h''(t) = 1 + O(t) \qquad\textup{for } t\to 0^+;
\end{equation}

\item  $f$ is convex and the matrix $ \Gradp^2 f(p)$ is uniformly positive definite in compact subsets, i.e. fulfills the following local ellipticity condition
\begin{equation}\label{ellipticity}
\forall M>0 \ \exists \lambda=\lambda(M)>0\quad:\quad \langle \Gradp^2 f(p)\xi, \xi\rangle \geq \lambda |\xi|^2 \quad\forall\;  |p|\leq M \mbox{ and } \forall \xi\in \mathbb{R}^n.
 \end{equation}
\end{itemize}

This class of problems obviously contains the linear case of the Dirichlet energy $h(t) = t^2/2$ and the geometric case of minimal surfaces $h(t) = \sqrt{1+t^2}-1$, where the constant $1$ is clearly irrelevant in the minimization problem.
The existence and uniqueness of a solution $u$ in the class $\mathcal{A}_g$ under suitable growth conditions on $f$ can be established following classical results by Giaquinta-Modica and Giusti \cite{GiaMod, Giu}. The solution to \eqref{min} can be characterized as the weak solution to the system
\begin{subequations}\label{Prob}\begin{align}
	\diver (\Gradp f\left(\nabla u \right))&=0  \mbox{ in } B^{+}_{1}, \label{equation}\\
	u \left( \Gradp f (\nabla u) \cdot e_{n} \right)&=0  \mbox{ on } B'_{1}, \label{bordo}\\
	-\Gradp f (\nabla u) \cdot e_{n} &\geq 0 \mbox{ on } B'_{1},\label{bordo2}\\
	u&\geq 0 \mbox{ on } B'_{1},\\
	u&= g \mbox{ on } (\partial B_1)^+,
 \end{align}\end{subequations}
where $e_n$ denotes the last vector of the standard basis of $\R^n$.

\subsection{Previous results}
The minimization problem \eqref{min} with the Dirichlet integral $\int |\nabla u|^2$ stems from the pioneering works of Signorini in elasticity theory \cite{Sig1}. Since then, lower dimensional obstacle problems naturally appeared in several fields, with applications for instance in fluid mechanics when describing osmosis through semi-permeable membranes, as well as in the boundary heat control problem, and many other contexts (see, e.~g., the book by Duvaut-Lions \cite{DuLi} and the survey papers \cite{FeRe, Ros2, Sal} for an extensive bibliography and further applications).

\medskip

Given for granted the existence of Lipschitz solutions, which in great generality have been shown to exist in the works by Giaquinta-Modica \cite{GiaMod} and Giusti \cite{Giu}, the main questions concern two aspects of the solutions $u$:
\begin{itemize}
\item the regularity of $u$ up to $B_1'$, where the boundary conditions are not prescribed but are implicitly determined by the solution itself through its free boundary,

\item the regularity of the free boundary, that is the subset of $B_1'$ where the solution ceases to saturate the one-side constraint and takes the natural boundary conditions.
\end{itemize}

\medskip

Most of the results in the literature deals with the linear case $\Gradp f(p) = p$, for which both the optimal regularity of the solution $u$ and the structure of the free boundary are known. As far as the first aspect is concerned, we recall the results by Lewy \cite{Lewy} and Richardson \cite{Richardson} for $n=2$, Caffarelli \cite{Caf}, Kinderlehrer \cite{Kin}, Ural\cprime tseva \cite{Ural} for the $C^{1,\alpha}$ regularity, for some $\alpha >0$; whereas the optimal regularity of the solutions is due to Athanasopoulos-Caffarelli \cite{AtCa}, where the authors show that the solution to \eqref{min} with $\Gradp f(p) = p$ are $C_{\textup{loc}}^{1,\frac12}$ in $B_1^+\cup B_1'$, and in general cannot be more regular.

\medskip

On the other hand, the nonlinear case is much less understood, despite its obvious relevance in the applications. One of the first occurrences of such instance can be found in a paper by Nitsche \cite{Nit} for the case of the minimal surface operator, for which the existence of Lipschitz solutions has been  first proven by Giusti \cite{Giusti2}. Since then, only partial results have been established for the nonlinear case, which turned out to be considerably more complicated than the linear counterpart. The first important contributions are given by Frehse \cite{Fre, Fre2}, where under very general assumptions the author shows the continuity of the first derivatives along the directions parallel to $B'_1$ in any dimension $n\geq 3$, and the continuity of the normal derivative in the case $n=2$ for the Lipschitz solutions to  nonlinear elliptic variational inequalities.
More refined results, regarding both the optimal regularity of the solutions and their free boundaries, have been found by Athanasopoulos \cite{Ath} for the minimal surface case in dimension $n=2$, and only recently more precise results for some nonlinear cases in general dimension appeared. We recall the case of minimal surfaces with flat obstacles treated by Focardi and the third named author \cite{FoSp20}, where the optimal H\"{o}lder continuity of the derivatives is proven, together with a study of the structure and the properties of the free boundary; as well as the study of a general class of nonlinear variational inequalities by Di Fazio and the third named author in \cite{DF-Sp-22}, where for the boundary obstacle problem the non-optimal $C^{1,\alpha}$ regularity in $B_1^+ \cup B'_1$ is shown; and the case of nonlinear thin obstacle problems considered
in \cite{DSV} in relation with the boundary Bernoulli problem in dimension $n=2$, where the optimal $C^{1,\frac12}$ regularity and the structure of the free boundary are derived with the use of complex analysis techniques.

\medskip

Our understanding is even more limited concerning the properties of the free boundary of the solutions: indeed, if for the linear problem many details on the free boundary have been investigated (structure of regular points \cite{CSS,FoSp16}, singular points \cite{GP}, non-regular and non-singular points \cite{FoSp18, FS22, AFS23}, only to mention few of the most relevant literature), for the nonlinear case much less is known, left aside the aforementioned results in the very special case of two dimensions and for the minimal surface case considered in \cite{FoSp20}.

\subsection{Main result: Optimal $C^{1, \frac{1}2}$ regularity}

The aim of this paper is to establish the optimal $C^{1, \frac{1}2}$-regularity for the general class of problems \eqref{min}.
The starting point is the $C^{1,\alpha}$-regularity established in \cite{DF-Sp-22}.

\begin{theorem}\label{t.main}
Let $u$ be the solution to the boundary obstacle problem \eqref{min} or equivalently \eqref{Prob}. Then,
$$u\in C_{\rm loc}^{1, 1/2}(B_1^+ \cup B_1'). $$
\end{theorem}

The idea of the proof is very much influenced by the work on stationary two-valued graphs done by Simon-Wickramasekera \cite{SiWi}. This connection has been already used in the case of the minimal surface operator in \cite{FoSp20}, whose $C^{1, \alpha}$ solutions fulfilled the hypotheses (suitably formulated) of \cite{SiWi}. Here instead we show how the proof of \cite{SiWi} can be adapted in order to cover the more general cases under examination.

We proceed in several steps. First, we introduce a new frequency function tailored to the nonlinear operator $\Gradp f(\nabla u)$ that is a modification of the famous \emph{Almgren frequency function} used for Dirichlet energy minimizers \cite{Alm} (see also \cite{DLSp-memo}). In order to show its (almost) monotonicity, we need to use the $C^{1,\alpha}$-regularity established in \cite{DF-Sp-22} and to exploit a dichotomy argument first used in \cite{SiWi}. The monotonicity formula for this modified frequency allows then to study the limit of the blow-up sequence and leads to the optimal $C^{1, \frac{1}2}$-regularity as in the classical case of quadratic energies. 

\medskip

A byproduct of our main result Theorem \ref{t.main} is that the results for the free boundary of the solutions to thin obstacle problems for a linear differential operator with H\"{o}lder continuous coefficients apply (see \cite{GaPeSV}), and moreover, if one has that $h(t) = \frac{t^2}{2} + O(t^4)$ (as in the minimal surface case, for instance), then the optimal regularity in Theorem \ref{t.main} implies that the solution $u$ solves a thin obstacle problem for a linear operator with Lipschitz coefficients, and more refined results on the structure of the singular set can be derived (see \cite{GaSm,FoSp20, AFS23}).
Moreover, we expect that the techniques developed in the present paper, and in particular the seek for a form of the frequency functions tailored to nonlinear problems, can be useful in related questions in the geometric calculus of variations (see, e.g., the general nonlinear energies for multiple valued functions \cite{DL-F-S}).

The more general nonlinear thin obstacle problems, starting from the question raised by Nitsche on minimal surfaces with non-constant unilateral thin constraints, remain still open and more refined techniques needs to be developed.

%

\vskip10pt
\textsc{Acknowledgment:}  The authors have been partially funded by the ERC-STG HiCoS ``Higher Co-dimension Singularities: Minimal Surfaces and the Thin Obstacle Problem'' Grant n. 759229.
 \vskip10pt

\section{Frequency function}\label{Freq}

\subsection{Schauder estimates}
The starting point of our analysis is the $C_{\textup{loc}}^{1,\alpha}(B^+_1 \cup B'_1)$-regularity of the solutions $u$ of \eqref{min} established in \cite{DF-Sp-22}.
For our aims, we need to establish the following Schauder estimates (not explicitly derived in \cite{DF-Sp-22}): there exist constants $\bar\varrho,k>0$ depending on the dimension $n$ and the Lipschitz constant of $u$ such that
\begin{equation}\label{Sch}
\|u\|_{L^{\infty}(B_{\varrho/2}^+)}+\varrho\|\nabla u\|_{L^{\infty}(B_{\varrho/2}^+)} + \vr^{1+\alpha} [\nabla u]_{\alpha, B^+_{\vr/2}}\leq k \left(\varrho^{-n} \int_{B_{\varrho}^+} u^{2}\right)^{\frac{1}{2}} \mbox{ for any } \varrho\in(0, \bar{\vr}),
\end{equation}
where given a radius $r>0$ the symbol $B_r^+$ denotes $B_r \cap \{x_n>0\}$, $B_r$ denotes the ball of $\R^n$ with radius $r$  centered at the origin, and $$[\nabla u]_{\alpha, B_\vr^+}:=\sup_{x\neq y\in B_\vr^+}\frac{|\nabla u(x) - \nabla u(y)|}{|x-y|^\alpha}$$ is the usual H\"{o}lder seminorm. 
Moreover, we show that $u\in W_{\textup{loc}}^{2,2}(B_1^+)$  and fulfills 
\begin{equation}\label{Der2}
\|\nabla^2u\|_{L^2(B_{\vr/2}^+)}\leq \beta \vr^{-2}\|u\|_{L^2(B_{\vr}^+)} \mbox{ for any } \vr\in (0, \bar{\vr}),
\end{equation}
where the positive constant $\beta$ depends on the dimension $n$ and on the Lipschitz constant of $u$.


\subsection{Frequency function}
We use the notation $\Lambda_u$ for the contact set:
\[
\Lambda_u:= \big\{u=0\big\}\cap B_1';
\]
and $\Gamma_u$ for the free boundary
\[
\Gamma_u:= \Big\{x\in B_1': B'_r(x) \cap \Lambda_u\neq \emptyset\quad \textup{and}\quad B'_r(x) \setminus \Lambda_u\neq \emptyset \quad \forall\;r>0\Big\}.
\]

In the sequel of this section we will always assume that the origin is a point of the free boundary:
\[
0\in \Gamma_u.
\]
Recall that by the $C^{1,\alpha}$-regularity this implies that
\[
u(0) = 0=|\nabla u (0)|.
\]
For any $\vr\in (0,1)$ we introduce the following nonlinear modification of the standard frequency function of $u$:
\begin{gather} \label{N}
N(\varrho):=\frac{D(\varrho)}{H(\varrho)},\\
\textup{where }\quad
D(\varrho):=\varrho^{2-n}\int_{B^{+}_{\varrho}}\Gradp f(\nabla u) \cdot\nabla u \d x,\qquad
H(\varrho):=\varrho^{1-n}\int_{(\pB{\varrho})^{+}} u^{2} \d \sigma_x.\notag
\end{gather}
For convenience we set the notation
\[
\|u\|_\vr = \left(\vr^{-n} \intS{u^2}\right)^\frac{1}{2}.
\]
The main result about the frequency is that, under the assumption of a \emph{doubling condition} for $u$,  a quasi-monotonicity formula  holds.

\begin{proposition}
\label{Prop-Mon} 
Let $u\in C_{\textup{loc}}^{1,\alpha}(B^+_1 \cup B'_1)$ be the solution to the boundary obstacle problem \eqref{Prob}, and assume that $0\in \Gamma_u$ and \eqref{Sch} and \eqref{Der2} hold with constants $\bar\varrho, k, \beta$.
If there is $\gamma>1$ and $\sigma\in (0, 1/2]$ such that
\begin{equation}\label{doubling0}
\|u\|_\sigma >0, \qquad \|u\|_\vr\leq \gamma \|u\|_{\frac{\vr}{2}} \quad\mbox{ for any } \vr \in (0, \sigma], 
\end{equation}
then there exists $\vr_0=\vr_0(\bar{\vr}, \sigma) \in (0,1)$ such that for all $\vr\in (0,\vr_0)$ the frequency function $N(\vr)$ is well-defined, bounded and 
\begin{equation} \frac{\d}{\d\vr}[{\rm exp}(\alpha^{-1} C \vr^\alpha)N(\vr)]\geq 0\quad \mbox{ for any } \vr \in (0,\vr_0),
\end{equation}
 where $C=C(n, \alpha, \gamma, \beta, k, [\nabla u]_\alpha)>0$.
\end{proposition}

\begin{remark}
Let us observe that the doubling condition \eqref{doubling0} can be written as
\begin{equation}\label{Doub}
\intS{u^2} \leq 2^n \gamma^2 \int_{B_{\vr/2}^+}{u^2\dx} \quad \mbox{ for any } \vr\in (0, \sigma],
\end{equation}
and in particular $\|u\|_\sigma>0$ implies that
\begin{equation}\label{well-def}
\|u\|_\vr>0  \mbox{ for any } \vr\in (0, \sigma].
\end{equation}
\end{remark}

In order to prove Proposition \ref{Prop-Mon} we need some technical results. 

\begin{remark}
Note that if $h$ satisfies \eqref{e.h''}, then  there exists $\bar t>0$ and functions $\omega_1$, $\omega_2$ such that
\begin{align}\label{der h1}
h''(t) = 1+ \omega_2(t),\qquad h'(t)=t\left(1 +\omega_{1}(t) \right),
\end{align}
with the functions $\omega_1$ and $\omega_2$ satisfying
\begin{align}
|\omega_2(t)| + |\omega_{1}(t)| + t |\omega_1'(t)|
\leq Ct\quad \forall\; t\in (0 ,\bar{t}),\qquad
\omega_2(t) = \omega_1(t) + t \omega_1'(t).
\end{align}
Moreover, formula \eqref{formula-f} yields
	\begin{align}\label{Df}
	\Gradp f(p) &=h'(|p|)\frac{p}{|p|}= (1 +\omega_{1}(|p|))\, p\\
	\label{D2f}
	\Gradp^{2}f(p)&=\frac{h'(|p|)}{|p|}\Id+\left( h''(|p|)-\frac{h'(|p|)}{|p|} \right) \frac{p\otimes p }{|p|^{2}}= (1 +\omega_{1}(|p|))\Id+  \omega_{1}'(|p|)  \frac{p\otimes p }{|p|}.
	\end{align}
	
Therefore, straightforward computations imply the following estimates:
\begin{align}\label{fD1}
|\nabla u -\Gradp f(\nabla u)|&\leq C |\nabla u|^{2},\\
|\Delta u -\diver(\Gradp f(\nabla u))|&\leq C |\nabla u||\nabla^{2}u|,\label{fD2}
\end{align}
where the constant $C$ is independent of $u$. 
\end{remark}

The following version of the Poincar\'e inequality is needed [see also \cite[Section $6$]{SiWi}].
\begin{lemma}
Let $u\in C^{1,\alpha}(B_1^+)$ with $u(0)=0$ fulfilling the Schauder estimates \eqref{Sch}. Then, 
\begin{equation}\label{Poincare}
\int_{B^+_{\varrho}} u^{2} \d x \leq C\varrho^{2}\int_{B^+_{\varrho}} |\nabla u|^{2} \d x \quad\mbox{ for any } \vr\in (0, \bar{\vr}),
\end{equation}
with $C>0$ depending on the dimension $n$, the constants in the Schauder estimates and in the doubling condition.
\end{lemma}

\bProof
Let $\vr\in (0, \bar{\vr})$, $\theta \in (0,1/2)$ and $\lambda:= \fint_{B_{\theta\varrho}^+}u \dx$, then the appropriate version of the Poincar\'e inequality (see \cite[Lemmas 7.12 and 7.16]{GiTr}) implies
\begin{equation}
\intS{|u-\lambda|^2}\leq C \vr^2 \intS{|\nabla u|^2} \quad \mbox{ with } C=C(n, \theta)>0.
\end{equation}
Therefore,
\begin{equation}\label{Poincare2}
\intS{u^2} \leq 2 \intS{\big(|u-\lambda|^{2}+\lambda^{2}\big)} \leq 2C \vr^2 \intS{|\nabla u|^2} +\omega_n \vr^n\lambda^2.
\end{equation}
Since there exists $y \in B_{\theta\varrho}^+$ with  $u(y)=\lambda$, and since $u(0)=0$, it follows that
\begin{equation}
\lambda\leq \theta\varrho \sup_{B^{+}_{\theta\vr}}|\nabla u|.
\end{equation}
Employing also the $L^\infty$-estimate for $\nabla u$ in \eqref{Sch}, from \eqref{Poincare2} we get that		
\begin{equation}
\intS{u^{2}}  \leq 
C \varrho^{2} \intS{ |\nabla u|^{2}} + k^2 \omega_n \theta^{2} \intS{ u^{2}}, \end{equation}
where $k$ is the constant in \eqref{Sch}.
This in turn gives \eqref{Poincare} choosing $\theta$ sufficiently small depending on $k$.
\qed
\\

In order to establish the boundedness of the frequency function we need the following auxiliary lemma. 
\begin{lemma}\label{lemma-ineq}
Let $u\in C_{\textup{loc}}^{1,\alpha}(B_1^+\cup B'_1)$ be the solution of the boundary obstacle problem \eqref{min} fulfilling the Schauder estimates \eqref{Sch} and \eqref{Der2} with constants $\bar \rho, k, \beta$, the doubling condition \eqref{Doub} with constants $\gamma, \sigma$, and $0\in \Gamma(u)$.
Then, it holds
\begin{align}
&\int_{(\pB{\varrho})^+} u^{2} \d\sigma_x \leq \frac{C_1}{\varrho}\int_{B_{\varrho}^+}u^{2}\d x \ \quad \forall\vr<\frac12\min \{\bar \vr, \sigma \}, \label{ineq2}
\end{align}
with $ C_1=C(n,k,\gamma)>0$ and there exists $\tilde \varrho \in (0, \bar\varrho)$ such that
\begin{align}\label{ineq1}
&\int_{B_{\varrho}^+} u^{2} \d x \leq \varrho\int_{(\pB{\varrho})^+}u^{2}\d\sigma_x +C_2\varrho^{1+ \alpha}\int_{B_{\varrho}^+} |u||\nabla u|\d x\ \quad \forall \vr\in (0, \tilde{\vr}),
\end{align}
with $C_2=C_2([\nabla u]_\alpha)>0$.
\end{lemma}

\bProof
We start proving \eqref{ineq2}. Employing the Schauder estimates \eqref{Sch} and the doubling condition \eqref{Doub} we have that for any $0<s<\vr<\frac12 \min \{\bar \vr, \sigma \}$
\[ \begin{split}
\frac{\d}{\d s}\left(\int_{(\pB{s})^+}u^{2} \d \sigma_x\right)&
=\frac{n-1}{s}\int_{(\pB{s})^+} u^{2}\d\sigma_ x + 2 \int\limits_{(\pB{s})^+}u\nabla u\cdot \eta \d \sigma_x\\
		&\leq \frac{(n-1)\, k^2}{s\, \varrho^n} \,|(\pB{s})^+|  \int_{B_{2\varrho}^+} u^{2}\d x + \frac{2 k^2}{\varrho^{n+1}}\,|(\pB{s})^+| \int_{B_{2\varrho}^+} u^{2} \d x\\
		&\leq (n-1)\,n\,2^n\,\omega_n \, s^{n-2}  \frac{k^2 \gamma}{\varrho^{n}} \int_{B_{\varrho}^+} u^{2}\d x + n\,\omega_n\,2^n\,s^{n-1} \frac{k^2\gamma}{\varrho^{n+1}} \int_{B_{\varrho}^+} u^{2} \d x\\
&\leq C(n,k,\gamma) \frac{1}{\varrho^{2}} \int_{B_{\varrho}^+} u^{2}\d x,
\end{split} \]
where $\eta(x) = \frac{x}{|x|}$ is the outer normal vector and $C(n,k,\gamma)>0$ is a constant.
Integrating with respect to $s\in (0,\varrho)$ with $\varrho <\frac12\min \{\bar \vr, \sigma \}$, \eqref{ineq2} follows.

Turning to \eqref{ineq1}, we notice that
\begin{align*}
\frac{\d}{\d s}\left(\int_{(\pB{s})^+}u^{2} \d \sigma_x\right)
&=(n-1)s^{n-2}H(s) + 2 \int_{(\pB{s})^+}u\nabla u \cdot \eta \d \sigma_x\\
&\geq 2 \int_{(\pB{s})^+}u \nabla u \cdot \eta \d \sigma_x.
\end{align*}
Since $u$ fulfills \eqref{Prob},  using the divergence theorem and \eqref{Df}, we deduce that
\begin{equation}\label{int uDu}
\begin{split}
\int_{(\pB{s})^+}u \nabla u \cdot \eta \d \sigma_x &= 
\int_{(\pB{s})^+}u \big(\nabla u - \nabla_p f(\nabla u)\big) \cdot \eta \d \sigma_x
+ \int_{(\pB{s})^+}u \nabla_p f(\nabla u) \cdot \eta \d \sigma_x\\
& = -\int_{(\pB{s})^+}u \,\omega_1(|\nabla u|) \nabla u \cdot \eta \d \sigma_x\\ 
&+ \int_{B_{s}^+}\diver\big(u \nabla_p f(\nabla u)\big) \d x
+\int_{B_s'}u \nabla_p f(\nabla u) \cdot e_n \d x'.
\end{split}
\end{equation}
Recall that $|\omega_1(t)|\leq C t$ for $t\leq \bar t$. Therefore, considering that $\nabla u(0)=0$, we find a radius $\tilde \varrho$ such that
\[
|\nabla u(x)|\leq \bar t\qquad\forall\;x\in B_{\tilde\varrho}^+.
\]
This implies that
\[
\left\vert
\int_{(\pB{s})^+}u \,\omega_1(|\nabla u|) \nabla u \cdot \eta \d \sigma_x
\right\vert \leq C \int_{(\pB{s})^+}|u|\,|\nabla u|^2 \d \sigma_x \qquad \forall\; s<\tilde \varrho.
\]
As far as the third integral in \eqref{int uDu} is concerned, we notice that by the Signorini boundary condition \eqref{Prob} we infer that
\[
\int_{B_s'}u \nabla_p f(\nabla u) \cdot e_n \d x' = 0.
\]
Finally, the second integral in \eqref{int uDu} is positive for small enough radii: indeed, 
\begin{align*}
\int_{B_{s}^+}\diver\big(u \nabla_p f(\nabla u)\big) \d x & = 
\int_{B_{s}^+}\nabla u \cdot \nabla_p f(\nabla u) \d x +
\int_{B_{s}^+}u\,\diver\big(\nabla_p f(\nabla u)\big) \d x\\
&= \int_{B_{s}^+}\nabla u \cdot \nabla_p f(\nabla u) \d x\\
&=\int_{B_{s}^+}|\nabla u|^2 \big(1+\omega_1(\nabla u)\big) \d x\geq 0,
\end{align*}
if $\tilde \rho$ is small enough to ensure that $|\nabla u(x)|\leq C^{-1}$ for $x\in B_{s}^+$.

We can then estimate as follows
\begin{equation*}\begin{split}
\frac{\d}{\d s}\left(\int_{(\pB{s})^+}u^{2} \d \sigma_x\right) &
\geq -C \int_{(\pB{s})^+} |u| |\nabla u|^{2} \d \sigma_x \geq -C s^{\alpha} \int_{(\pB{s})^+} |u| |\nabla u|\d\sigma_x 
\end{split}\end{equation*}
where we employed the H\"{o}lder continuity of $\nabla u$. Integrating the last inequality for $s\in (\tau, \varrho)$, with $ \vr<\tilde \vr $, we obtain
\[ \int_{\tau}^{\varrho}	\frac{\d}{\d s}\left(\int_{(\pB{s})^+}u^{2} \d\sigma_ x\right) \d s\geq -C \varrho^{\alpha} \int_{\tau}^{\varrho}	 \int_{(\pB{s})^+} |u| |\nabla u|\d\sigma_x \d s \]
which gives
\[ \int_{(\pB{\varrho})^+}u^{2} \d\sigma_x +C \varrho^{\alpha}  \int_{B_{\varrho}^+} |u| |\nabla u|\d x\geq \int_{(\pB{\tau})^+}u^{2} \d \sigma_x. \]
Integrating again with $\tau \in (0,\varrho)$
$$\varrho \int_{(\pB{\varrho})^+}u^{2} \d x +C \varrho^{\alpha+1} \int_{B_{\varrho}^+} |u| |\nabla u| \d x  \geq  \int_{B_{\varrho}^+}u^{2} \d x$$
which is \eqref{ineq1}.
\qed

\subsection{Boundedness of the frequency.}

If $0\in \Gamma_u$, then $u$ is a nonzero solution of the thin obstacle problem and this implies that  $\int_{(\pB{\varrho})^+}u^{2}\d \sigma_x>0$ for every $\varrho>0$.

By virtue of \eqref{ineq2} and the Poincar\'e inequality \eqref{Poincare} it follows that for any $\vr\in (0, \frac12\min \{\sigma,\bar{\vr}\})$
\begin{equation}
 \int_{(\pB{\varrho})^+}u^{2}\d \sigma_x \leq \frac{C}{\varrho}\int_{B_{\varrho}^+}u^{2}\d x 
		\leq  C\varrho \int_{B_{\varrho}^+}|\nabla u|^{2}\d x. 
 \end{equation}
On the other hand, by the Schauder estimates \eqref{Sch} and the doubling condition \eqref{Doub} we get 
\begin{equation}\label{Nbelow}
 \varrho^2\int_{B_{\varrho}^+}|\nabla u|^{2} \d x \leq k^2 \left( \frac{1}{\varrho^{n}} \int_{B_{2\varrho}^+}u^{2}\d x\right)|B_{\varrho}^+|
 \leq  C\int_{B_{\varrho}^+}u^{2}\d x,
 \end{equation}
 where $ C=C(n,k,\gamma) $.
Next \eqref{ineq1} and the Cauchy-Schwartz inequality imply that
 \begin{equation*}\begin{split}
		\int_{B_{\varrho}^+} u^{2} \d x &\leq \vr\int_{(\pB{\varrho})^+}u^{2}\d\sigma_x +C\varrho^{1+ \alpha }\int_{B_{\varrho}^+} |u||\nabla u|\d x\\
		& \leq\vr \int_{(\pB{\varrho})^+}u^{2}\d\sigma_x +C\varepsilon\varrho^2\int_{B_{\varrho}^+} |\nabla u|^{2}\d x+\frac{C\varrho^{2\alpha}}{\varepsilon}\int_{B_{\varrho}^+} |u|^{2}\d x
	\end{split} \end{equation*}
for any $\varepsilon>0$. Then there exists $\vr_0(\varepsilon)\in (0,\tilde \vr)$ such that
\[ \left( 1- \frac{C\varrho^{2\alpha}}{\varepsilon} \right) \int_{B_{\varrho}^+} u^{2} \d x  \leq \vr \int_{(\pB{\varrho})^+}u^{2}\d\sigma_ x +C\varepsilon\varrho^2\int_{B_{\varrho}^+} |\nabla u|^{2}\d x \]
for any $\vr\in (0, \vr_0(\varepsilon))$ (with positive left-hand side). Going back into \eqref{Nbelow} it yields
\[ \varrho^2\int_{B_{\varrho}^+}|\nabla u|^{2} \d x \leq C\vr \int_{(\pB{\varrho})^+}u^{2}\d \sigma_x +C\varepsilon\varrho^2\int_{B_{\varrho}^+} |\nabla u|^{2}\d x \]
and thus choosing $\varepsilon>0$ sufficiently small we conclude
\begin{equation}\label{DminH}
\varrho\int_{B_{\varrho}^+}|\nabla u|^{2} \d x \leq C \int_{(\pB{\varrho})^+}u^{2}\d x. 
\end{equation} 
Therefore using that 
$$ 0< c_1\leq 1 + \omega_1(|\nabla u|) \leq c_2  \qquad\forall\;x\in B_{\vr_0}^+ \mbox{ with } \vr_0 \mbox{ sufficiently small},$$
we can conclude that, if $0\in \Gamma_u$ and the doubling condition \eqref{Doub} holds, then there exist positive constants $C_1, C_2,\varrho_0$ (depending on $n$, $k$, $\gamma$, $[\nabla u]_\alpha$) such that
\begin{equation}\label{Nlimitata}
C_1\leq N(\vr)\leq C_2 \mbox{ for any } \vr\in (0, \vr_0].
\end{equation}

\subsection{Quasi-monotonicity of the frequency.}
Here we prove Proposition \ref{Prop-Mon}.

Deriving \eqref{N} we obtain
\begin{equation}
N'(\varrho)= \frac{D'(\vr)H(\vr)-D(\vr)H'(\vr)}{H(\vr)^2},
\end{equation}
with
\begin{align}
D'(\vr)&=(2-n)\vr^{1-n}\intS{\Gradp f(\nabla u)\cdot\nabla u} + \vr^{2-n}\intB{\Gradp f(\nabla u)\cdot\nabla u}, \label{Dprimo}\\
H'(\vr)&=2\vr^{1-n}\intB{u u_\vr}\label{H1},
\end{align}
where $u_\vr= \nabla u \cdot x/\vr $ denotes the radial derivative. Now we can rewrite $N'(\vr)$ as following 
\begin{equation}\label{CS}\begin{split}\
N'(\varrho)&=
 \frac{1}{(H(\vr))^2}\left[\left( 2 \vr^{2-n} \intB{|u_\vr|^2}\right) H(\vr) - \left(\vr^{2-n} \intB{u u_\vr} \right)H'(\vr) \right]\\
 &+ \frac{1}{(H(\vr))^2}\left[ \left(D'(\vr) - 2 \vr^{2-n} \intB{|u_\vr|^2}\right)H(\vr) \right.\\
&\qquad\qquad\qquad \left.- \left(D(\vr) -  \vr^{2-n} \intB{u u_\vr} \right) H'(\vr)\right].
\end{split}\end{equation}
Then the Cauchy-Schwartz inequality implies that 
\begin{equation}\label{N1}\begin{split}
N'(\vr)&\geq \frac{1}{(H(\vr))^2}\left[ \left(\!D'(\vr) - 2 \vr^{2-n} \intB{|u_\vr|^2}\!\right)H(\vr)\right.\\
&\qquad\qquad\qquad\left. - \left(\!D(\vr) -  \vr^{2-n} \intB{u u_\vr}\! \right) H'(\vr)\right].
\end{split}\end{equation}
We now estimate each term in the right-hand side of \eqref{N1}, beginning with 
$$ E_1:= \left| D'(\vr) - 2 \vr^{2-n} \intB{|u_\vr|^2}\!  \right|.  $$

%
%

Using the divergence theorem first and $x\cdot e_n=0$ for every $x\in B_\vr'$, we deduce that
\begin{equation}\label{divergenza}\begin{split}
\vr\intB{\Gradp f(\nabla u)\cdot \nabla u}&=\intS{\diver(\Gradp f(\nabla u) \cdot\nabla u\, x )}+\int_{B_\vr'} \Gradp f(\nabla u)\cdot \nabla u \, x\cdot e_n \d x\\
&=n \intS{\Gradp f(\nabla u) \cdot\nabla u} + \intS{(D^2u \Gradp f(\nabla u))\cdot x} \\
&\qquad+  \intS{ (\partial_{p_l}\Gradp f(\nabla u)\cdot \nabla u) \  (\partial_l\nabla u\cdot x)}.
\end{split}\end{equation}
First let us consider the second integral on the right-hand side and integrate by parts (with respect to the $j$-th variable in $D^2 u=\partial_i\partial_j u$)
\begin{equation*}\begin{split}
\intS{(D^2 u \Gradp f(\nabla u)) \cdot x}&
= \intB{\nabla u\cdot x\, \Gradp f(\nabla u)\cdot \nu}-\int_{B_\vr'} \nabla u\cdot x\, \, \Gradp f(\nabla u)\cdot e_n \d x\\
&\qquad-\intS{\Gradp f(\nabla u) \cdot \nabla u}-\intS{(\nabla u \cdot x)\, \big(D^2_{p}f(\nabla u) : D^2 u\big)}\\
&= \intB{\nabla u\cdot x\, \Gradp f(\nabla u)\cdot \nu}-\intS{\Gradp f(\nabla u) \cdot \nabla u}\\
&\qquad-\intS{(\nabla u \cdot x)\,\big( D^2_{p}f(\nabla u): D^2 u\big)},
\end{split}
\end{equation*}
where we employed the Signorini boundary condition \eqref{bordo} on $B'_\vr$ in the boundary integral at the right-hand side.
Now, using \eqref{Df} and \eqref{D2f}, we get
\begin{equation}\label{Abb}
\begin{split}
\intS{(D^2 u &\Gradp f(\nabla u)) \cdot x}\\
& =\intB{(1+\omega_{1}(|\nabla u|)) \nabla u\cdot x\, \nabla u\cdot \nu}-\intS{\Gradp f(\nabla u) \cdot\nabla u}\\&-\intS{\nabla u \cdot x\,(1+ \omega_1(|\nabla u|))\Delta u } -\intS{\nabla u \cdot x\, \omega_{1} '(|\nabla u|)\,  \frac{\nabla u\otimes \nabla u}{|\nabla u|}: D^2 u}.
\end{split}
\end{equation}
Next the third term at the right-hand side of \eqref{divergenza}:
\begin{equation}
\begin{split}
&\intS{ (\partial_{p_l}\Gradp f(\nabla u)\cdot \nabla u) \  (\partial_l\nabla u\cdot x)}= \intS{\partial^2_{p_l p_j} f(\nabla u) \partial_j u\, \partial^2_{l i}u \, x_i}\\
	&=\intS{(1+\omega_{1}(|\nabla u|)) (D^2 u \nabla u) \cdot x} + \intS{\omega_{1}'(|\nabla u|) |\nabla u| (D^2u \nabla u)\cdot x }=:A + B.
\end{split}
\end{equation}
We integrate by parts in $A$ with respect to the $ j-$th variable in  $D^2 u=\partial_i\partial_j u$, and use \eqref{bordo} on $B'_\vr$:
\begin{equation}\label{Abb2}\begin{split}
&A= \intS{(1+\omega_{1}(|\nabla u|)) \partial_{ij}u \, \partial_j u\,  x_i} =\intB{(1+\omega_{1}(|\nabla u|)) \,  \nabla u\cdot \nu\,\nabla u\cdot x}\\
&-\intS{\Gradp f(\nabla u) \cdot\nabla u}-\intS{\nabla u \cdot x\,(1+ \omega_1(|\nabla u|))\Delta u }\\
&\qquad -\intS{\nabla u \cdot x\, \omega_{1} '(|\nabla u|)\,  \frac{\nabla u\otimes \nabla u}{|\nabla u|}: D^2 u}.
\end{split}\end{equation}

Collecting \eqref{divergenza}-\eqref{Abb2} we realize that

\begin{equation}
	\begin{split}
\vr\intB{\Gradp f(\nabla u) \cdot\nabla u}
&=(n-2) \intS{\Gradp f(\nabla u) \cdot\nabla u}\\
&\qquad +2 \intB{(1+\omega_{1}(|\nabla u|)) \,  \nabla u\cdot \nu\,\nabla u\cdot x}\\
&\qquad-2\intS{\nabla u \cdot x\,(1+ \omega_1(|\nabla u|)\Delta u }\\
& \qquad-2\intS{\nabla u \cdot x\, \omega_{1} '(|\nabla u|)\,  \frac{\nabla u\otimes \nabla u}{|\nabla u|}\cdot D^2 u} + B.
\end{split}
\end{equation}

Therefore it results that
\begin{equation}\begin{split}
 D'(\vr) - 2 \vr^{2-n} \intB{\left |\nabla u \cdot \frac{x}{\vr}\right|^2}&
 = 2 \varrho^{1-n}\intB{\omega_{1}(|\nabla u|) \nabla u\cdot x\,  \nabla u\cdot \nu}\\
 &\qquad-2\varrho^{1-n}\intS{\nabla u \cdot x\,(1+ \omega_1(|\nabla u|)\Delta u } \\ 
 &\qquad-2 \varrho^{1-n}\intS{\nabla u \cdot x\, \omega_{1} '(|\nabla u|)\,  \frac{\nabla u\otimes \nabla u}{|\nabla u|}\cdot D^2 u} \\
 &\qquad+  \varrho^{1-n} \intS{\omega_{1}'(|\nabla u|) |\nabla u| (D^2u \nabla u)\cdot x }\\
 &=:\sum_{i=1}^4 A_i,
\end{split}\end{equation}

hence
\begin{equation}\begin{split}
E_1:&=\left| D'(\vr) - 2 \vr^{2-n} \intB{|u_r|^2} \right| \leq \sum_{i=1}^4 |A_i|.
\end{split}\end{equation}
Let us discuss each term $|A_i|$. It holds:
\begin{equation*}
|A_1|\leq C\vr^{2-n+\alpha} \|\nabla u\|_\infty^2 |\partial B_\vr| \stackrel{\eqref{Sch}}{\leq} C\vr^{-1-n+\alpha}  \intSS{ u^{2}} \stackrel{\eqref{Poincare}}{\leq} C\vr^{1- n + \alpha}  \intS{ |\nabla u|^{2}},
\end{equation*}

\begin{equation*}
\begin{split}
|A_2|&\leq C\vr^{1-n} (1+C\vr^{\alpha}) \intS{|\nabla u| \, |x||\Delta u-\diver \Gradp f(\nabla u)| }\\ 
&\stackrel{\eqref{fD2}}{\leq}C \vr^{2-n+ \alpha}\intS{|D^2u| |\nabla u|}\\
&\leq C \vr^{2-n+\alpha}\left(\intS{|\nabla u|^2}\right)^{\frac{1}{2}} \left(\intS{|D^2u|^2}\right)^{\frac{1}{2}} \\
&\stackrel{\eqref{Der2}}{\leq}
C \vr^{-n+\alpha}\left(\intS{|\nabla u|^2}\right)^{\frac{1}{2}}\left(\intSS{|u|^2}\right)^{\frac{1}{2}}\\
&  \stackrel{\eqref{Doub},\;\eqref{Poincare}}{\leq} C \vr^{1-n+\alpha}\intS{|\nabla u|^2}. 
\end{split}
\end{equation*}
Similarly,
\begin{align*}
\begin{split}
|A_3|+ |A_4|&\leq \vr^{1-n}\intS{ |\omega_{1}'(|\nabla u|) | |D^2u|  |\nabla u|^{2} |x|} \leq C\vr^{2-n +\alpha} \intS{|D^{2}u| |\nabla u|}\\
& \leq C \vr^{1-n+\alpha}\intS{|\nabla u|^2}. 
\end{split}\end{align*}

Therefore,
\begin{equation}\label{E1}
E_1\leq C \vr^{1-n+\alpha}\intS{|\nabla u|^2}\leq C \vr^{\alpha-1}D(\vr).
\end{equation}

Now, let us consider
\begin{equation}\label{E2} E_2:=  \left|D(\vr) -  \vr^{2-n} \intB{u u_\vr}\right|.\end{equation}

Integrating by parts and employing \eqref{bordo} we get
\begin{equation*}
	\begin{split}
	\intS{\Gradp f(\nabla u)\cdot \nabla u}&=\intBp{u\, \Gradp f(\nabla u)\cdot \nu}-\intS{ u \,\partial_{j} (\partial_{p_{j}}f(\nabla u)) }\\
	&=\intB{u\, (1+ \omega_1(|\nabla u|))\nabla u \cdot \nu } - \intS{u\, D^2f(\nabla u): D^2u}\\
	&= \intB{u\, (1+ \omega_1(|\nabla u|))\nabla u \cdot \nu }  -\intS{(1+ \omega_1(|\nabla u|))\, u \,\Delta u} \\&-\intS{u \, \omega_{1}'(|\nabla u|)\, \frac{\nabla u\otimes\nabla u}{|\nabla u|}: D^2u}.
	\end{split}
\end{equation*}
Thus
\begin{equation}\begin{array}{l}\displaystyle\vspace{6pt}
	E_{2} =\left| \vr^{2-n} \intB{u\, \omega_1(|\nabla u|)\nabla u \cdot \nu } - \vr^{2-n}\intS{(1+ \omega_1(|\nabla u|))\, u \,\Delta u} \right.\\
	\displaystyle\vspace{6pt}\hfill\left.-\vr^{2-n}\intS{u \, \omega_{1}'(|\nabla u|)\, \frac{\nabla u\otimes\nabla u}{|\nabla u|}: D^2u}\right|,
\end{array}	\end{equation}
then using that $u\in C^{1,\alpha}(B_\vr^+)$
\begin{equation}\label{}
	\begin{split}	 
	E_2\leq  \vr^{3-n +\alpha} \|\nabla u\|_\infty^2 |\partial B_\vr| + C \vr^{2-n}(1+C\vr^{\alpha}) \intS{| u| |\Delta u -\diver D_p f(\nabla u)|} &\\
	+C\vr^{2-n}\intS{|u| |D^{2}u| |\nabla u|},&
	\end{split}
\end{equation}
and employing the estimate \eqref{fD2}, the H\"{o}lder inequality, the Schauder estimates \eqref{Der2} together with the doubling assumption and the Poincar\'{e} inequality \eqref{Poincare} we conclude that 

\begin{equation}\label{E22}
E_2\leq C\vr^{2-n+\alpha}\intS{|\nabla u|^2} \leq C\, \vr^{\alpha} \,D(\vr).
\end{equation}
Moreover, by virtue of \eqref{H1} and \eqref{E2} we can deduce that
\begin{equation}\label{H11}
|H'(\vr)|\leq \left|H'(\vr)- \frac{2}{\varrho} D(\vr)\right| + \frac{2}{\varrho}|D(\vr)|\leq C\vr^{-1}(E_2 +|D(\vr)|).
\end{equation}
Putting together \eqref{Nlimitata},\eqref{E1}, \eqref{E22}, \eqref{H11} into \eqref{N1}  we conclude that
\begin{equation}
N'(\vr)\geq - C\vr^{\alpha-1}N(\vr).
\end{equation}
Whence we infer that
\begin{equation}
\frac{\d}{\d\vr}({\rm exp}(\alpha^{-1}C\vr^\alpha)N(\vr))\geq 0,
\end{equation}
thus concluding the proof of Proposition \ref{Prop-Mon}.

\section{Blowup sequence}\label{s3}
Here we prove that the rescalings at a free boundary point of the solutions to the nonlinear obstacle problems \eqref{min} converge to the solutions of the classical Signorini problem for the Dirichlet energy.

\begin{proposition}\label{p.blowup}
Let $u\in C_{\textup{loc}}^{1,\alpha}(B^+_1 \cup B'_1)$ be the solution to the boundary obstacle problem \eqref{Prob}, $0\in \Gamma_u$, 
and let \eqref{Sch} and \eqref{Der2} hold with constants $\bar\varrho, k, \beta$.
Assume that there is $\gamma>1$ and $\sigma\in (0, 1/2]$ such that 
\begin{equation}\label{doubling}
\|u\|_\sigma >0, \quad \|u\|_\vr\leq \gamma \|u\|_{\frac{\vr}{2}} \mbox{ for any } \vr \in (0, \sigma].
\end{equation}
Then, for every sequence of positive real numbers $\vr_{j}$ such that $\vr_{j} \searrow 0,$  there is a subsequence $ \vr_{j'} $ such that 
$$ \frac{\vr_{j'}^{\frac{n}{2}}u(\vr_{j'}x)}{\|u\|_{L^{2}(B^{+}_{\vr_{j'}})}} \rightarrow \varphi(x) \quad\mbox{ in } C^1 \mbox{ locally on } \{x_n\geq 0\},$$
where $\varphi\in C^{1,\frac12}(\{x_n>0\})$ is a homogeneous solution to the Signorini problem:
\[
\int_{(\partial B_1)^+}|\nabla \varphi|^2\, \d x \leq
\int_{(\partial B_1)^+}|\nabla \psi|^2\, \d x \qquad \forall\; \psi\vert_{(\partial B_1)^+} = \varphi\vert_{(\partial B_1)^+}, \; \psi\vert_{B_1'}\geq 0.
\]
Moreover,
\begin{equation}
e^{C \vr^{\alpha}}N_{u}(\vr)\geq N_{u}(0)=N_{\varphi}(0)\geq \frac{3}{2}\quad \mbox{ for any } \vr \in (0, \vr_0), \end{equation}
and
\begin{equation}\label{under}
 \|u\|_{\vr} \leq C\left(\frac{\vr}{\sigma}\right)^{\frac{3}{2}}  \|u\|_{\sigma} \quad\mbox{ for any } \vr \in (0,\sigma],
 \end{equation}
 with $C=C(n, k, \gamma, [\nabla u]_\alpha)>0$ and $\vr_0$ is the same as in Proposition \ref{Prop-Mon}.
\end{proposition}

\bProof
Let $\vr_j$ be a sequence of real numbers such that $\vr_j \searrow 0$ and let $u_j$ be the sequence of functions defined as $u_j(x)=\lambda_j u(\vr_j x)$ with
$$\lambda_j:=\frac{\vr_{j}^{\frac{n}{2}}}{\|u\|_{L^{2}(B^{+}_{\vr_{j}})}}$$ is chosen such that $\|u_j\|_{L^2(B_1^+)}=1$ for any $j\in \mathbb{N}$. 
By virtue of \eqref{Sch} and \eqref{doubling} (see also \eqref{Doub}) we have that $\|u_{j}\|_{C^{1,\alpha}(B_{R}^+)}\leq C(n,R, \gamma) $ for all $R>0$ and for $j$ sufficiently  large. Thus there is a subsequence $u_{\vr_{j'}} $ that converges in $C^{1}$ locally on $\{x_n\geq 0\}$ to a function $ \varphi=\varphi(x) \in C^{1, \alpha}$. With an abuse of notation we keep denoting the subsequence with $u_{\vr_{j}} $.

Let us show that $\varphi$ is solution of the Signorini problem. For any $\psi$ suitable test function with compact support in $B_1^+$ it holds that 
\begin{equation}\label{uj}
\int_{B_1^+}{\nabla u_j\cdot\nabla \psi\dx} +
 \int_{B_1^+} \omega_1\left(\textstyle{\frac{|\nabla u_j|}{\lambda_j \vr_j}}\right) \, \nabla u_j \cdot \nabla \psi\dx=0 \end{equation}
because $u$ solves \eqref{Prob}.
Since 
$$\frac{|\nabla u_j(x)|}{\lambda_j \vr_j} = |\nabla u(\varrho_j x)| \leq C\varrho_j^\alpha\qquad \forall\; x\in B_1^+,
$$
then 
\begin{equation}\label{limit-rhoj} 
\omega_1\left(\textstyle{\frac{|\nabla u_j|}{\lambda_j \vr_j}}\right)\to 0\quad 
\mbox{ as } j \to +\infty.\end{equation}
Thus taking the limit in \eqref{uj} and employing the convergence in $C^{1}$ of $u_j$ to $\varphi$ we get
$$\int_{B_1^+}{\nabla\varphi \cdot \nabla \psi\dx} =0,$$
that means that $\varphi$ is harmonic in $B_1^+$.
Next, let us discuss the limit of the conditions on $B_1'$.
Since $u$ solves \eqref{Prob}, then $u_j$ solves 
\begin{equation}\label{thinobst-limit}
u_j \Gradp f\left (\frac{\nabla u_j}{\lambda_j\vr_j}\right) \cdot e_n =0 
\mbox{ on } B_1'.\end{equation}
Note that
\begin{align*}
\Gradp f\left (\frac{\nabla u_j}{\lambda_j\vr_j}\right) &= h'\left (\frac{|\nabla u_j|}{\lambda_j\vr_j}\right)\frac{\nabla u_j}{|\nabla u_j|}
= \left(1+\omega_1\left (\frac{|\nabla u_j|}{\lambda_j\vr_j}\right)\right) \frac{\nabla u_j}{\lambda_j\vr_j}.
\end{align*}
Therefore, 
$$ \lambda_j\vr_j \Gradp f\left (\frac{\nabla u_j}{\lambda_j\vr_j}\right) \cdot e_n=
\left(1+\omega_1\left (\frac{|\nabla u_j|}{\lambda_j\vr_j}\right)\right) \nabla u_j \cdot e_n \to \nabla \varphi \cdot e_n,$$
and hence we get
\begin{equation}
\varphi\, \nabla \varphi\cdot e_n =0
\qquad\textup{and}\qquad \nabla \varphi\cdot e_n\leq 0\quad \mbox{ on } B_1',
\end{equation}
concluding that $\varphi$ is a solution to the classical Signorini problem.
For any radius $\sigma>0$, we denote $N_j(\sigma)$ the frequency for the sequence $u_j$ and it holds that
$$N_j(\sigma)= \frac{(\sigma\vr_j)^{2-n} \int_{B_{\vr_j \sigma}^+} |\nabla u (y)|^2 \d y}{(\sigma\vr_j)^{1-n} \int_{\partial B_{\vr_j \sigma}^+} |u (y)|^2 \d \sigma_y}=N_{u}(\varrho_j\sigma),$$
whence 
\begin{equation}
N_j(\sigma) \to \lim_{j\to +\infty} N_u(\varrho_j\sigma) = N_u(0).
\end{equation}
On the other hand the convergence in $C^{1}$ ensures that
\begin{equation}
N_j(\sigma) \to N_\varphi(\sigma) \mbox{ for any radius } \sigma.
\end{equation}
Combining the last two, it means that  $N_\varphi(\sigma)$ is constant and by standard results we infer that $\varphi$ is homogeneous of degree $N_u(0)\geq 3/2$ (see \cite{AC,ACS}).

In order to establish \eqref{under} we start noticing that 
as in \eqref{H11} and \eqref{E22} we know that 
\[  H'(\vr) \geq C(1-\vr^{\alpha})\frac{1}{\vr}D(\vr), \]
hence  employing  the monotonicity of $N(\vr)$
\[ (\ln H(\vr))'=\frac{H'(\vr)}{H(\vr)}\geq C\frac{1}{\vr}N(\vr) \geq C\frac{1}{\vr}N(0). \]
Then, integrating onto the interval $(\vr , \sigma)$, it yields 
\[ \ln \left(\frac{H(\sigma)}{H(\vr)}\right)\geq C\ln \left( \frac{\sigma}{\vr}\right)^{\frac{3}{2}} \]
which means that
\begin{equation}\label{H/H}
\left( \frac{H(\sigma)}{H(\vr)}\right)^{\frac{1}{2}}\geq C \left( \frac{\sigma}{\vr}\right)^{\frac{3}{2}}. 
 \end{equation}
 
 By virtue of \eqref{ineq2} it holds 
\begin{equation}\label{H(s)}
 H(\sigma)=\sigma^{1-n} \int_{\pB{\sigma}^+} u^{2} \d \sigma_x\leq C\sigma^{-n} \int_{B_{\sigma}^+} u^{2} \d x,
 \end{equation}
while on the other side from \eqref{DminH} and \eqref{Poincare} we have
\begin{equation}\label{H(r)}
 H(\vr)\geq C \vr^{-n} \int_{B_{\vr}^+} u^{2} \d x.
 \end{equation}
Therefore \eqref{H/H} implies the conclusion
\[  \left(\vr^{-n} \int_{B_{\vr}^+} u^{2} \d x\right)^{\frac{1}{2}}\leq C\left(\frac{\vr}{\sigma}\right)^{\frac{3}{2}}\left(\sigma^{-n} \int_{B_{\sigma}^+} u^{2} \d x\right)^{\frac{1}{2}}.  \]
\qed

\section{Proof of the optimal $C^{1, \frac{1}{2}}$ regularity}\label{Blowupseq}

Finally in this section we show the optimal regularity $u\in C_{\textup{loc}}^{1, 1/2}(B_1^+\cup B'_1)$.
We fix the notation
\[
\|u\|_{\varrho, x_0} := \left(\varrho^{-n} \int_{B_{\varrho}^+(x_0)} u^2\,\d x\right)^{\frac12}, \qquad x_0\in \Gamma_u.
\]
We begin with the following preliminary lemma (cp. \cite{SiWi}).

\begin{lemma}\label{lem6.1} Let $u$ be solution to \eqref{min} and let $\delta\in (0,\frac12)$.
Then, there exists $\vr_{0}\in (0, {1}/{2}],$
such that for every $x_0\in \Gamma_u\cap B_{\frac12}$  and every $\vr\in (0,\vr_{0}]$ the following implication holds:
$$
\|u\|_{\frac{\vr}{2},x_0}\geq 2^{\frac{3}{2}+\delta} \|u\|_{\frac{\vr}{4},x_0} \quad\Longrightarrow \quad\|u\|_{\sigma,x_0}\geq \left(\frac{2\sigma}{\vr}\right)^{\frac{3}{2} + \delta} \|u\|_{\frac{\vr}{2},x_0} \qquad
\forall\; \sigma\in\left [\textstyle{\frac{3}{4}}\vr,\vr\right] .$$
\end{lemma}

\bProof We argue by contradiction. Let us suppose that for any $l\in \mathbb{N}\setminus\{0,1\}$ there are $ \vr_{0}=\frac{1}{l},$ $\vr_{l}\leq \frac{1}{l}$, $\sigma_l\in [\frac34\vr_l, \vr_l]$ and solutions $u_l=u$ to \eqref{min} with $x_l\in \Gamma_{u_l}$  such that the assertion does not hold.
Then, introducing the rescaling $$\tilde{u}_{l}(x):=\|u_{l}\|_{2\vr_{l},x_l}^{-1}u_{l}(2\vr_{l} x+x_l),$$ by the estimates \eqref{Sch} we have that, up to a not relabeled subsequence, the functions $\tilde{u_{l}}$ converge in $C^{1}$ to $\varphi$ on $ B^+_{\sigma} $ with $ \sigma=\lim_l\frac{\sigma_l}{2\vr_l}\in [\frac{3}{8},\frac12]$ that is a solution to the Signorini problem in $B_1$.
Let us now consider the function $ f(\vr)=\vr^{-3/2 - \delta}\|\varphi\|_{\vr} $, with $ 0<\vr\leq \sigma $ and $\delta>0$. From the assumption, we have that the function $f$ attains its maximum in $ [\frac{1}{8}, \sigma] $ in an interior point $ \tilde{\vr} \in (\frac{1}{8},\sigma )$: indeed, 
$$
\|u_l\|_{\frac{\vr_l}{2}}\geq 2^{\frac{3}{2}+\delta} \|u_l\|_{\frac{\vr_l}{4}}
\quad \Longrightarrow \quad
 \left({\frac14}\right)^{-\frac32 - \delta }\|\tilde u_l\|_{\frac{1}{4}}\geq  \left({\frac18}\right)^{-\frac32 - \delta}
\|\tilde u_l\|_{\frac{1}{8}} \quad \Longrightarrow \quad
f\left({\frac14}\right)\geq f\left({\frac18}\right),
$$
and analogously
\begin{align*}
\|u_l\|_{\sigma_l}< \left(\frac{2\sigma_l}{\vr_l}\right)^{\frac{3}{2} + \delta} \|u_l\|_{\frac{\vr_l}{2}}
&\quad \Longrightarrow \quad
\left(\frac{\sigma_l}{2\vr_l}\right)^{-\frac32 -\delta}\|\tilde u_l\|_{\frac{\sigma_l}{2\vr_l}}< \left({\frac14}\right)^{-\frac32 -\delta} \|\tilde u_l\|_{\frac{1}{4}}\\
&\quad \Longrightarrow \quad
f\left(\sigma\right)\leq f\left({\frac14}\right).
\end{align*}
Rather than $f'(\vr)$ let us compute $(f^2(\vr))'$:
\begin{align*}
\left( \vr^{-3-2\delta - n}\intS{\varphi^{2}}\right)'&= \left( \vr^{-3-2\delta}\int_{B_1^+}\varphi^{2}(\varrho x)\, \d x\right)'\allowbreak
 \\
 &=(-3 -2\delta)\vr^{-3 -1-2\delta-n}\intS{\varphi^{2}}+ 2\vr^{-3-2\delta-n}\intS{\varphi_{r}\varphi}\allowdisplaybreaks\\
 & =(-3 -2\delta) \vr^{-3-1-2\delta-n}\int_{0}^{\vr}\int_{\pB{t}^+} \varphi^{2} \d \sigma_x \d t\\
 &\quad + 2\vr^{-3-2\delta-n} \int_{0}^{\vr}\int_{\pB{t}^+} \varphi_{r}\varphi \d \sigma_x \d t\allowdisplaybreaks\\
 & =(-3 -2\delta)\vr^{-3-2\delta-2}\int_{0}^{\vr}H_\varphi(t) \d t+ 2\vr^{-3-2\delta-2} \int_{0}^{\vr}D_\varphi(t) \d t\allowdisplaybreaks\\
 &=2 \vr^{-3-2\delta-2}\int_{0}^{\vr}H_\varphi(t)\left(N_\varphi(t)-\frac{3}{2}-\delta\right) \d t.
\end{align*}
Since $\varphi$ is a solution to the Signorini problem, the monotonicity of $ N $ and the fact that $N_\varphi\geq \frac32$ imply that $ f $ is monotone increasing.
Moreover, considering that under the contradiction assumption, $f$ must have an interior maximum,  we infer that $f$ needs to be constant in the interval $I=[\frac{1}{8}, \sigma]$. Therefore, we conclude that actually $ N_{\varphi}(\vr)=3/2 + \delta$ for $ \vr \in I $.
This implies that $ \varphi$ is a homogeneous solution of degree $3/2 + \delta<2$, thus contradicting the fact that there are no solution with frequency $\lambda \in (\frac32, 2)$ (see \cite{ACS}).\qed

\begin{remark}
Let us notice that for $ \sigma=\frac{\vr}{2} \in (0,\frac{\vr_0}{2}) $ in Lemma \ref{lem6.1} we have that
\[ \|u\|_{\sigma,x_0} \geq 2^{\frac{3}{2}+\delta} \|u\|_{\frac{\sigma}{2},x_0} \quad\Longrightarrow \quad\|u\|_{2\sigma,x_0}\geq 2^{\frac{3}{2} + \delta} \|u\|_{\sigma,x_0} \qquad. \]
It follows that, if $\|u\|_{\sigma,x_0}\geq 2^{\frac{3}{2} + \delta} \|u\|_{\sigma/2,x_0}$ then $\|u\|_{2^j \sigma,x_0} \geq 2^{\frac{3}{2} + \delta} \|u\|_{2^{j-1}\sigma,x_0}$ for every $j=1,2,\dots$ such that $2^{j}\sigma\leq \vr_0$, which implies
\begin{equation}
\|u\|_{\vr,x_0}\leq C \vr^{\frac{3}{2} + \delta } \|u\|_{\vr_0,x_0} \quad\mbox{ for any } \vr\in (\sigma/2, \vr_0/2],
\end{equation}
and, by changing the constant $C$, also
\begin{equation}\label{sopra}
\|u\|_{\vr,x_0}\leq C\vr^{\frac{3}{2}} \|u\|_{L^2(B_1)}  \quad \mbox{ for any } \vr\in (\sigma/2, 1/2].
\end{equation}
\end{remark}

\begin{proposition}\label{p.decay}
Let $u\in C_{\textup{loc}}^{1,\alpha}(B_1^+\cup B_1')$ be the solution of \eqref{min} and $x_0\in \Gamma_u\cap B_{\frac12}$. Then,
\begin{equation}\label{e.decay}
\sup_{B^+_{\vr/2}(x_0)} |u| + \vr\sup_{B_{\vr/2}^+(x_0)} |\nabla u| \leq  C\vr^{\frac{3}{2}}\|u\|_{L^{2}(B_1^+)} \quad \forall\;\vr\in (0,{1}/{2}),
\end{equation}
where $C=C(n, \alpha, \beta, k)>0$.
\end{proposition}

\bProof
Assume without loss of generality that $x_0=0\in \Gamma_u$.
Let us introduce
\[  \sigma := \inf\left\lbrace \frac{1}{2}, \lbrace \vr \in (0,{1}/{2}] : \|u\|_{\vr}\geq 2^{\frac{3}{2} + \delta}  \|u\|_{\frac{\vr}{2}} \rbrace \right\rbrace, \]
where $ \delta \in (0, \frac12) $.
From Lemma \ref{lem6.1} and \eqref{sopra},  we have
\begin{equation}\label{6.25}
\|u\|_{\vr}  \leq C\vr^{\frac{3}{2}}\|u\|_{L^{2}(B_1^+)}\quad  \mbox{ for any } \vr\in (\sigma/2, 1/2].
\end{equation}
On the other hand, if  $\vr \in (0, \sigma] $, then it holds
\[ \|u\|_{\vr}< 2^{\frac{3}{2} + \delta}  \|u\|_{\frac{\vr}{2}}  \]
and thus by Proposition \ref{p.blowup} and \eqref{under}
\begin{equation}\label{6.26}
\|u\|_{\vr}\leq C\left(\frac{\vr}{\sigma}\right)^{\frac{3}{2}}  \|u\|_{\sigma} \mbox{ for any } \vr \in (0, \sigma].
\end{equation}
Combining \eqref{6.25} and \eqref{6.26} we deduce that there exists a constant $C>0$ such that
\[ \|u\|_{\vr}  \leq C\vr^{\frac{3}{2}}\|u\|_{L^{2}(B_1)}\quad \mbox{ for any } \vr\in (0,{1}/{2}). \]
By \eqref{Sch}, we have finally
\[ \sup_{B_{\vr/2}^+} |u| + \vr\sup_{B_{\vr/2}^+} |\nabla u| \leq C \|u\|_{\vr}\leq  C\vr^{\frac{3}{2}}\|u\|_{L^{2}(B_1^+)} \quad \mbox{ for any } \vr\in (0,{1}/{2}).  \]
\qed

\begin{theorem}\label{theo6.5}
Every solution $u$ of \eqref{min} belongs to $C^{1,\frac12}_{\textup{loc}}(B_1^+\cup B_1')$.
\end{theorem}

\bProof
Recall that by \cite{DF-Sp-22} every solution of \eqref{min} is $C_{\textup{loc}}^{1, \alpha}$. The local $C^{1,\frac12}$ regularity of $u$ is now an easy consequence of Proposition \ref{p.decay}. We give the proof for readers' convenience.
We need to estimate the oscillation of the gradient.
We distinguish between three cases.\\
\textsc{Case 1.} Let $z, y\in B^+_{1/2}$ be such that
\[
\textup{dist}(z, \Gamma_u)\leq 10^{-1}|z-y|\qquad\textup{and}\qquad \textup{dist}(y, \Gamma_u)\leq 10^{-1}|z-y|.
\]
Then, let $p_z, p_y\in \Gamma_u$ be such that
\[
\textup{dist}(z, \Gamma_u)=|z-p_z|, \qquad \textup{dist}(y, \Gamma_u)=|y-p_y|,
\]
and note that
\begin{align*}
|\nabla u (z) - \nabla u(y)|&\leq |\nabla u (z) - \nabla u(p_z)|+|\nabla u (p_z) - \nabla u(p_y)|+|\nabla u (p_y) - \nabla u(y)|\\
& = |\nabla u (z) - \nabla u(p_z)|+|\nabla u (p_y) - \nabla u(y)|\\
& \stackrel{\eqref{e.decay}}{\leq} C |z-p_z|^{\frac12} + C |y-p_y|^{\frac12}\leq C|z-y|^\frac{1}{2}.
\end{align*}

\textsc{Case 2.} Let $z,y\in B^+_{1/2}$ be such that
\[
10^{-1}|z-y|<\varrho_z:=\textup{dist}(z, \Gamma_u), \qquad B_{\varrho_z}(z)\cap \Lambda_u = \emptyset.
\]
We consider the symmetrization matrix $S\in\R^{n\times n}$, $Sx := (x',-x_n)$ for $x=(x',x_n)$, and set
\[
v(x) :=
\begin{cases}
u (Sx) & x_n\leq 0,\\
u(x) & x_n>0,
\end{cases}
\qquad x\in B_{\varrho_z}(z).
\]
It is easy to verify that $v$ satisfies
\[
\textup{div} \big(\nabla_pf(\nabla v)\big) = 0 \qquad \textup{in } B_{\varrho_z}(z).
\]
Indeed, if $x_n<0$, then
\[
\nabla v (x) = S \nabla u(Sx), \qquad \nabla^2v(x) = S\nabla^2u(Sx) S.
\]
Therefore, 
\begin{align*}
\textup{div} \big(\nabla_pf(\nabla v(x))\big) & = \nabla^2_pf\big(\nabla v(x)\big) : D^2 v(x) \\
&= \left[\big(1+\omega_1(|\nabla v(x)|)\big) I + \omega'_1(|\nabla v(x)|)\frac{\nabla v(x) \otimes \nabla v(x)}{|\nabla v(x)|}\right] :D^2v(x)\\
&= \left[\big(1+\omega_1(|\nabla u(Sx)|)\big) I + \omega'_1(|\nabla u(Sx)|)\frac{S\nabla u(Sx) \otimes S\nabla u(Sx)}{|\nabla u(Sx)|}\right] :SD^2u(Sx)S\\
&= \big(1+\omega_1(|\nabla u(Sx)|)\big) I:D^2u(Sx)	\\	
&\qquad + \omega'_1(|\nabla u(Sx)|)\frac{S\nabla u(Sx) \otimes S\nabla u(Sx)}{|\nabla u(Sx)|}:D^2u(Sx)\\
&= \textup{div} \big(\nabla_pf(\nabla u(Sx))\big) = 0.
\end{align*}
Moreover, for every $\varphi\in C_c^1(B_{\varrho_z}(z))$ we have that 
\begin{align*}
\int_{B_{\varrho_z}(z)}\nabla_pf(\nabla v(x)) \cdot \nabla \varphi(x) \,\d x & =
\int_{B_{\varrho_z}^+(z)}\nabla_pf(\nabla u(x)) \cdot \nabla \varphi(x) \,\d x\\
&\qquad+\int_{B_{\varrho_z}(z)\cap \{x_n<0\}}\nabla_pf(\nabla v(x)) \cdot \nabla \varphi(x) \,\d x\\
& = - \int_{B_{\varrho_z}'(z)} \nabla_p f(\nabla u(x',0))\cdot e_n \, \varphi(x',0) \, \d x'\\
&\qquad + \int_{B_{\varrho_z}'(z)} \nabla_p f(\nabla v(x',0))\cdot e_n \, \varphi(x',0) \, \d x' = 0,
\end{align*}
where we use that $\nabla_p f(p)= (1+\omega_1(|p|))p$ leads to
\[
\nabla_p f(\nabla v(x',0))\cdot e_n 
\nabla_p f(\nabla u(x',0))\cdot e_n = 0 \qquad \forall\; (x',0)\in B_1'\setminus \Lambda_u.
\]
Therefore, since $v$ solves a uniformly elliptic equation in $B_{\varrho_z}(z)$, by standard estimates we get 
\begin{align*}
|\nabla u (z) - \nabla u(y)| &\leq \|D^2 v\|_{L^\infty(B_{\varrho_z/2}(z))}|\nabla v (z) - \nabla v(y)|\\
&\leq C\varrho_z^{-2} \|v\|_{L^\infty(B_{\varrho_z}(z))}|\nabla v (z) - \nabla v(y)|\\
&\leq C \varrho_z^{-1/2}\,|z-y|\leq C |z-y|^{1/2}.
\end{align*}

\textsc{Case 3.} Let $z,y\in B^+_{1/2}$ be such that
\[
10^{-1}|z-y|<\varrho_z:=\textup{dist}(z, \Gamma_u), \qquad B_{\varrho_z}(z)\cap B_1'\subset \Lambda_u.
\]
We consider the symmetrization matrix $S\in\R^{n\times n}$, $Sx := (x',-x_n)$ for $x=(x',x_n)$, and set
\[
v(x) :=
\begin{cases}
- u (Sx) & x_n\leq 0,\\
u(x) & x_n>0,
\end{cases}
\qquad x\in B_{\varrho_z}(z).
\]
We need to verify as before that $v$ solves the equation:
\[
\textup{div} \big(\nabla_pf(\nabla v)\big) = 0 \qquad \textup{in } B_{\varrho_z}(z).
\]
Indeed, if $x_n<0$, then as above the equation is pointwise verified and, for every $\varphi\in C_c^1(B_{\varrho_z}(z))$, we have that 
\begin{align*}
\int_{B_{\varrho_z}(z)}\nabla_pf(\nabla v(x)) \cdot \nabla \varphi(x) \,\d x
& = - \int_{B_{\varrho_z}'(z)} \nabla_p f(\nabla u(x',0))\cdot e_n \, \varphi(x'0) \, \d x'\\
&\qquad + \int_{B_{\varrho_z}'(z)} \nabla_p f(\nabla v(x',0))\cdot e_n \, \varphi(x'0) \, \d x' = 0,
\end{align*}
where now we use that
\begin{align*}
\nabla_p f(\nabla v(x',0))\cdot e_n 
&= (1+\omega_1(|\nabla v(x',0)|))\partial_n v(x',0) \\
& = (1+\omega_1(|\nabla u(x',0)|))\partial_n u(x',0)\\
& = 
\nabla_p f(\nabla u(x',0))\cdot e_n \qquad\qquad \forall\; (x',0)\in B_1'\cap \Lambda_u.
\end{align*}
We can then conclude exactly as in the previous case that 
\begin{align*}
|\nabla u (z) - \nabla u(y)| \leq C |z-y|^{1/2},
\end{align*}
thus finishing the proof.
\qed

\bibliographystyle{amsplain}

\providecommand{\bysame}{\leavevmode\hbox to3em{\hrulefill}\thinspace}
\providecommand{\MR}{\relax\ifhmode\unskip\space\fi MR }
\providecommand{\MRhref}[2]{%
  \href{http://www.ams.org/mathscinet-getitem?mr=#1}{#2}
}
\providecommand{\href}[2]{#2}

\end{document}